\documentclass[10pt]{article}

\usepackage{amsthm,amsmath}
\usepackage{amssymb}
\usepackage[square,sort,comma,numbers]{natbib}
\usepackage{graphicx}

\usepackage[utf8]{inputenc}
\newcommand{\beginsupplement}{%
        \setcounter{table}{0}
        \renewcommand{\thetable}{S\arabic{table}}%
        \setcounter{figure}{0}
        \renewcommand{\thefigure}{S\arabic{figure}}%
     }

\allowdisplaybreaks

\begin{document}

\title{Multicriteria global optimization for biocircuit design}

\author{Irene Otero-Muras\thanks{ireneotero@iim.csic.es} and Julio R. Banga\thanks{julio@iim.csic.es}  \\ BioProcess Engineering Group, \\IIM-CSIC (Spanish Council of Scientific Research)\\
  Eduardo Cabello 6, 36208, Vigo, Spain }

\maketitle
\begin{abstract}
 One of the challenges in Synthetic Biology is to design circuits with increasing levels of complexity. While circuits in Biology are complex and subject to natural tradeoffs, most synthetic circuits are simple in terms of the number of regulatory regions, and have been designed to meet a single design criterion.

In this contribution we introduce a multiobjective formulation for the design of biocircuits. We set up the basis for an advanced optimization tool for the modular and systematic design of biocircuits capable of handling high levels of complexity and multiple design criteria. Our methodology combines the efficiency of global Mixed Integer Nonlinear Programming solvers with multiobjective optimization techniques. Through a number of examples we show the capability of the method to generate non intuitive designs with a desired functionality setting up a priori the desired level of complexity. The presence of more than one competing objective provides  a realistic design setting where every design solution represents a trade-off between different criteria. The tool can be useful to explore and identify different design principles for synthetic gene circuits.
\end{abstract}

\section*{Introduction}
A hallmark of Synthetic Biology is, quoting Arkin, the ambition \emph{to formalize the process
of designing cellular systems in the way that traditional engineering disciplines have formalized design and manufacture, so that complex
behaviours can be achieved for practical ends} \cite{Arkin:2006}. In formalizing the design process, as it is the case in more traditional engineering disciplines,
mathematical modeling and optimization play a central role.

Over the past ten years, many advances have been achieved in the field, from the first bacterial toggle switches \cite{Gardner:2000} and  biological oscillators \cite{Elowitz:2000},
to the recent mammalian cell to cell communication devices \cite{Bacchus:2012}. In a so called first wave of Synthetic Biology basic elements and small biological
modules were successfully implemented and characterized. One of the challenges of the second wave in progress is the integration of modules to create circuits of increasing complexity \cite{Purnick:2009}. However, as reported by Purnick and Weiss \cite{Purnick:2009} the level of complexity achieved in synthetic circuits, measured by the number of regulatory regions, is relatively low.  While circuits in Biology are complex, subject to natural tradeoffs and playing multiple roles \cite{Szekely:2013}, most synthetic designs are simple and perform a single task.
Reported reasons for the current limited complexity in synthetic circuits include too simplistic engineering design principles \cite{Purnick:2009}, difficulty to independently control multiple cellular processes in parallel \cite{Rao:2012} and increasing problems to construct and test circuits as they get larger \cite{Kwok:2010}.
Efforts are necessary to overcome this difficulties and, quoting Lu \emph{et al.} \cite{Lu:2009},
\emph{advancing Synthetic Biology to the realm of higher-order networks with programmable functionality and real world applicability}. In parallel, new computational tools need to be developed to support these efforts \cite{Medema:2012}.

In this contribution, our goal is to set up the basis of an advanced optimization tool for the \emph{modular} and \emph{systematic} design of biocircuits capable of handling high levels of complexity and multiple design criteria.

\emph{Modular} design requires the previous definition of standardized functional objects and interfaces \cite{Canton:2007}.
From the foundations of Synthetic Biology, efforts have been held in order to characterize \emph{standard biological parts}, i.e. DNA sequences encoding a function that can be assembled with other standard parts.
The abstraction hierarchy proposed by Endy \cite{Endy:2005} classifies standard parts in three different layers: \emph{parts},
that are defined as sequences with basic biological functions (like for example DNA-binding proteins), \emph{devices}
that are combinations of parts with a particular function and \emph{systems} that are combinations of devices.
An emerging catalogue of standard parts is available at the registry supported by the BioBricks Foundation \cite{biobricks}.

\emph{Systematic} design relies on mathematical models describing the circuit dynamics.  In this regard, modular modeling tools are advancing to facilitate the
mathematical representation of biological parts and their combinations \cite{Marchisio:2013}, providing the description of the reactions taking place inside the different parts and the interfaces to connect them.
Inspired by the BioBrick registry of standard parts, Marchisio and Stelling \cite{Marchisio:2008} developed a formal modeling framework based on the ODE formalism that permits modular model composition, that has been recently extended for the modeling of
more complex eukariotic systems \cite{Marchisio:2013}. Some remarkable advances have been also achieved regarding synthetic biology computer aided design tools \cite{Huynh:2013}.

The systematic design of circuits combining components or parts from a list or library can be formulated as an optimization problem \cite{Huynh:2013,Dasika:2008,Rodrigo:2007}
where the circuit model structure is manipulable through decision variables, and the desired behaviour of the circuit is encoded in the objective function to optimize.
This results in Mixed Integer Nonlinear Problems  (MINLP) whose solution is challenging due to the simultaneous
presence of binary variables and constraints in form of ODE's.

Dasika and Marnas \cite{Dasika:2008} developed an optimization framework for the design of biocircuits,  based on the circuit modeling formulation by Hasty \cite{Hasty:2001} and a multistart local outer approximation method for the optimization. A number of design problems were successfully solved within this framework including a circuit with inducer specific response, a genetic decoder and a concentration band detector.

In this work, we advance the optimization-based design of biocircuits with two contributions:
increasing the computation efficiency in order to handle higher levels of complexity and introducing multiple criteria in the design. To this purpose, we first introduce
a set of global MINLP solvers that reduce drastically the computation time for the monoobjective design problem in comparison with other published methods. Then we formulate a general multiobjective optimization framework that combines the efficiency of the global MINLP solvers with the ability to tackle multiple design criteria.
The inducer specific response circuit design by Dasika and Maranas \cite{Dasika:2008} is used to illustrate the efficiency of the MINLP methods presented and further reformulated with additional design criteria to discuss the advantages of a multiobjective formulation in the design of genetic circuits.

\section*{Methods}
\subsection*{Global Stochastic MINLP solvers for Biocircuit Design}
Optimization based design of biocircuits requires the integration of tools for modular modeling, simulation and optimization.
As reported in the introduction, modular tools for modeling in Synthetic Biology are advancing fast as well as repositories of biological parts. Searching for a generic optimization framework, the methods presented next do not bound to a specific modeling tool, but accommodate to any ODE based modeling framework such that the circuit's model structure can be obtained from the starting list of parts by giving values to a set of integer variables.

The design problem consists of finding the best solution or solutions among the set of all possible alternatives according to a number of criteria. In this first part, we focus on problems with one unique design objective.

Under these assumptions, the design of biocircuits can be formulated as a Mixed Integer Nonlinear Programming Problem \cite{Dasika:2008,Rodrigo:2007},  where
the model structure can be  encoded by integer variables and the constraints are the dynamics of the system in form of ODE's. Tunable kinetic parameters are real decision variables in the optimization model.
For a complete formulation we refer to Ref. \cite{Dasika:2008}, where the single objective MINLP problem is formalized for a particular modeling framework \cite{Hasty:2001}.

Next, our focus is on the computational challenges of the resultant MINLP, since some features inherent to biological circuit models make it particularly difficult to solve.

In first instance, the dynamics of biocircuits are highly nonlinear, and the resultant optimization problem  is non convex and multi-modal. In this type of problems, local methods lead to suboptimal solutions (unless we start close to the global optimum).  A number of approaches have been proposed in previous works to find the global optimum in monoobjective biocircuit design. Dasika and Maranas \cite{Dasika:2008} implemented a multistart local outer approximation algorithm
where a convergence sequence of upper and lower bounds to the original problem is generated and a local optimum solution is identified at each iteration. 
In this way, a local deterministic search is performed from several points.
Rodrigo \emph{et al.} \cite{Rodrigo:2007} use a stochastic metaheuristic based on simulated annealing  \cite{Kirkpatrick:1983,Cerny:1985} and
Huynh \emph{et al.} \cite{Huynh:2012} apply a global deterministic optimization method to a linear approximation of the model around a steady state.

On the other hand, the design of gene circuits involves in general large search spaces that combine a high number of
integer variables with the presence of real variables. Our first goal is to provide global optimization methods that efficiently solve monoobjective design problems of medium/high complexity.
Global deterministic methods ensure convergence to the global optimum within a desired tolerance, but the computational burden
is in general very high for non convex systems with large search spaces. Therefore, we have decided to employ global stochastic methods, which offer no guarantee of convergence to the global minimum in a finite number of iterations but showed excellent results solving complex process optimization problems in reasonable computation time \cite{Egea:2010}.

In this work, we use three different global stochastic methods: mixed-integer tabu search (\texttt{MITS}) \cite{Exler:2008}, mixed-integer ant colony optimization (\texttt{ACOmi}) \cite{Schlueter:2009} and the enhanced scatter search \texttt{eSS} described in \cite{Egea:2010}.
The three methods are actually hybrid, since the stochastic global search is combined with the local mixed-integer sequential quadratic
programming (\texttt{MISQP}) developed by Exler \emph{et al.} \cite{Exler:2007}. These methods have been shown to be efficient metaheuristics in solving complex-process optimization problems from different fields, providing a good compromise between
diversification (exploration by global search) and intensification
(local search).

\texttt{MITS} uses a combinatorial component, based on Tabu
Search \cite{Glover:1997}, to guide the search into promising areas, where the local solver is activated to precisely approximate local
minima.  Exler \emph{et al.} \cite{Exler:2008} made use of \texttt{MITS} to solve complex integrated design problems where other state of
the art solvers failed, including a wastewater plant for nitrogen removal and the well known Tennessee Eastman Process.

\texttt{ACOmi} extends ant colony optimization meta-heuristic \cite{Dorigo:2004} to handle mixed integer search domains. Schlueter \emph{et al.} \cite{Schlueter:2009} showed the efficiency of this method for a number of engineering benchmark problems with high levels of non-convexity.

\texttt{eSS} is an enhanced version of the scatter search by for mixed integer search domain.  Egea \emph{et al.} \cite{Egea:2010} proved the efficiency of the method for solving complex-process models through a set of engineering benchmarks, where \texttt{eSS} performed well even in cases in which standard local search methods failed to locate the global solution.

In this contribution, we evaluate the efficiency of these methods in the context of Synthetic Biology and in particular for the systematic design of genetic circuits.

For illustrative purposes we chose a representative design example from Ref. \cite{Dasika:2008}, with one single design objective.
Starting from a list of components, the goal is to build a circuit with a specific response upon stimulation by two different inducers.
There are eight different promoter elements (denoted by $P_1 \hdots P_8$):  $Plac_1$, $Plac_2$, $Plac_3$, $Plac_4$, $P_{\lambda}$, $Ptet_1$, $Ptet_2$ $Para$ and four transcripts (denoted by $R_1 \hdots R_4$):
$tetR$, $lacI$, $cI$, and $araC$. The inducers of interest are $IPTG$ and $aTc$. The dynamic model of the overall reaction network is constituted by a set of ordinary differential equations of the form:
\begin{equation}\label{eq:ex1balance}
\frac{z_j}{dt}(t)= V_j(t)  -{K_j}_{decay}  z_j (t)  ~~~\forall j
\end{equation}
where $V_j$ is the generation/consumption rate of $z_j$ due to the reactions and ${K_j}_{decay} z_j$ is the degradation rate. The rate expressions for the transcripts are known and they read:
\begin{equation}\label{eq:ex1rates}
V_j(t)=\sum_{i} Y_{ij} v_{ji}(t)
\end{equation}
where $v_{ji}$ is the rate of production of $R_j$ from $P_i$, and $Y_{ij}$ is a binary variable such that:
\begin{align*}
 Y_{ij}&=1 ~\textnormal{production of protein} ~R_j ~\textnormal{from promoter} ~P_i ~\textnormal{is turned on} \\
 Y_{ij}&=0 ~\textnormal{otherwise}
\end{align*}

The structure of the model is given by a superstructure $8 \times 4$ matrix $Y$  containing the $32$ binary variables of the model. We define the vector of binary variables $y$ as the vector obtained by converting the matrix $Y$ to
a vector by columns. The tunable parameters are contained in a vector of real variables denoted by $x$.

As mentioned, the goal is to achieve a specific response upon induction. Namely, the steady state level of $LacI$ must be high upon $aTc$ and low upon $IPTG$ induction whereas the steady state level of $tetR$ must be low upon $aTc$ and high upon $IPTG$ induction. This design goal is encoded in the following objective function $Z$ to be maximized:

\begin{equation}\label{eq:ex1Z}
Z=\left(\frac{[lacI]^{ss}_{aTC}-[lacI]^{ss}_{IPTG}}{[lacI]^{ss}_{IPTG}}+\frac{[tetR]^{ss}_{IPTG}-[tetR]^{ss}_{aTc}}{[tetR]^{ss}_{IPTG}}\right)/2
\end{equation}
where the maximum value $Z=1$ is achieved for $[lacI]^{ss}_{IPTG}=[tetR]^{ss}_{aTc}=0$.

The design problem is formulated as a MINLP where the decision variables are contained in the vectors $y$ and $x$,  the objective function to maximize is $Z$ in (\ref{eq:ex1Z}), subject to the system's dynamics (\ref{eq:ex1balance}).

The following constraint on the maximum number of active pairs is also imposed:
\begin{equation}
\sum_i \sum_j Y_{ij} \leq M_{max}.
\end{equation}
thus limiting the complexity of the circuit.

First we use the original formulation of the problem by Dasika and Maranas \cite{Dasika:2008}, with a maximum of two promoter-transcript pairs, and compare the performance of the methods with the published results. After we gradually increase the network complexity to
evaluate how the methods proposed scale with the increasing problem size. The results obtained are included in Results and Discussion section.

\subsection*{Multiobjective framework for automatic Biocircuit Design}
In traditional engineering disciplines design problems are often multicriteria, where a
number of design objectives are conflicting (typically production and cost) since we cannot increase one without decreasing the other.
Problems with multiple and conflicting design criteria do not have a unique optimal solution, but a trade-off front
between the competing objectives, also known as Pareto optimal front of solutions.

In biological systems, trade-offs between robustness, fragility, performance, and resource demands have been conjectured \cite{Szekely:2013,el2005optimal,sendin2006model,higuera2012multi, shoval2012evolutionary}. We know that living organisms allocate limited resources to various competing traits, and arising tradeoffs are central to evolutionary biology.  Furthermore molecular pathways have been shown in many cases to play diverse and complex roles.
However,  \emph{de novo} engineered circuits have designed to perform a single task and optimization based designs in Synthetic Biology have been formulated as problems with a single objective.

In this contribution we propose a multiobjective optimization framework for the design of biocircuits. In first instance, the design is formulated as a multicriteria optimization problem with a number of conflicting objectives and then a multiobjective optimization strategy is implemented to find the Pareto optimal set of solutions.

In order to mathematically define the multiobjective design problem, let first introduce the following vectors: $z \in \mathbb{R}^n$ is   the vector of state variables coding for the levels of all the species involved in the circuit;
$x \in \mathbb{R}^r$ is   the vector of continuous variables containing a set of tunable parameters; $y \in \mathbb{Z}^b$ is the vector of integer variables determining the circuit model structure;
$k \in \mathbb{R}^k$  is the vector of fixed parameters and  $J_i(\dot{z},z,x,y,k)$ for $i=1,\hdots,s$  is the set of conflicting objectives, where one subset of objectives encodes aspects related to the performance of the circuit and a second subset encodes aspects related to robustness and/or cost.

The design of a biocircuit can be formulated as finding a vector $x \in \mathbb{R}^{r}$ of continuous
variables and a vector $y \in \mathbb{Z}^{b}$ of integer variables which minimize the vector $J$ of $s$ objective functions:

\begin{subequations}\label{MOP}
\begin{equation}\label{MOP1}
\min_{x,y}    ~~  J_1 (\dot{z},z,x,y,k), J_2 (\dot{z},z,x,y,p), \hdots, J_s (\dot{z},z,x,y,k)
\end{equation}

subject to:
\begin{itemize}
\item[i)] {the circuit dynamics in the form of ODEs (or DAEs) with the state variables $z$ and additional parameters $k$:
\begin{equation}\label{MOP2}
f(\dot{z},z,x,y,k)=0,~~ z(t_0)=z_0,
\end{equation}

}
\item[ii)] {additional requirements in the form of equality and inequality constraints:
\begin{equation}\label{MOP3}
h(z,k,x,y)=0,
\end{equation}
\begin{equation}\label{MOP4}
g(z,k,x,y)\leq0,
\end{equation}
}
\item[iii)] {upper and lower bounds for the real and integer decision variables:
\begin{equation}\label{MOP5}
x_L~\leq~ x ~\leq~ x_U,
\end{equation}
\begin{equation}\label{MOP6}
y_L~\leq~ x ~\leq~ y_U.
\end{equation}
}
\end{itemize}
\end{subequations}
In order to evaluate the solutions of the multiobjective optimization problem, we need to introduce the notion of Pareto optimality \cite{Sendin:2010}.
Given two pairs $(x^*,y^*)$, $(x^{**},y^{**})$, we say that the vector $J(x^{*},y^{*})$ dominates $J(x^{**},y^{**})$ if $J(x^{*},y^{*}) \leq J(x^{**},y^{**})$ for all $i=1,\hdots,s$ with at least one strict inequality.
A feasible circuit defined by $(x^*,y^*)$ is a Pareto optimal solution of the multiobjective optimization problem if it is not dominated by other feasible circuits.
The set of all Pareto optimal solutions is known as Pareto front.

Computing the Pareto optimal set is a very challenging task in the context of complex biocircuit design. On the one hand, as indicated previoulsy, high complexity imply large search spaces, and on the other hand the expected Pareto front
is discrete and possibly non-convex, due to the high nonlinearity of the biocircuits dynamics and the existence of discrete decision variables.

There are a number of approaches to solve multiobjective optimization problems (MOPs) \cite{miettinen1999nonlinear}. Evolutionary approaches \cite{deb2001multi} allow to compute an approximation of the entire Pareto front in one single run, but require large population sizes
and consequently a high computational effort for the systems with the complexity we want to tackle.
Scalar approaches consist in transforming the MOP into one or more single objective problems, and include among others the well known \emph{weighted sum approach}, Normal Boundary Intersection (NBI) and $\epsilon$-constraint methods \cite{Sendin:2010}.

In the weighted sum approach, weights must be changed in order to generate different solutions in the Pareto front and the performance depends on the choice of the weighting coefficients, which is in general not straightforward. This method cannot find solutions in concave parts of the Pareto front.

NBI first builds a plane in the objective space which contains all convex combinations of the individual minima, denoted as convex hull of individual minima (CHIM) and then constructs normal lines to this plane.
The MOP is reformulated as to maximize the distance from a point on the CHIM along the normal through this point. When dealing with integer variables, there may not exist a feasible
solution on the selected normal to the CHIM, and therefore NBI at least in its original formulation has limited applicability for discrete Pareto fronts.

In the $\varepsilon$-constraint strategy \cite{Sendin:2010}, the MOP is reduced to a a number of MINLP, where each MINLP is obtained by minimising one of the objectives and converting the rest of criteria to inequality constraints.
Different solutions can be obtained by changing the upper bounds on the objectives not minimised. This methodology has two important advantages for the design of complex biocircuits:  the methodology works well for discrete and non-convex Pareto fronts and, in addition, it allows exploiting the MINLP solvers introduced in the previous section, that solve efficiently the resultant MINLPs at a reasonable cost.
Next we describe the $\varepsilon$-constraint strategy implemented in this work.  The proposed optimization process is composed of the following steps (for simplicity and without loss of generality we have considered two objective functions $J_1$ and $J_2$):

\begin{itemize}
\item[1.] {Search for the optima of each of the individual objectives:
 \[(x_1^*,y_1^*), ~~(x_2^*,y_2^*).\]}
\item[2.] {Compute the individual objective bounds as:
\[\underline{J_1}=J_1(x_1^*,y_1^*) , ~\overline{J_1}=J_1(x_2^*,y_2^*), \]
\[\underline{J_2}=J_2(x_2^*,y_2^*) , ~\overline{J_2}=J_2(x_1^*,y_1^*). \]}
\item[3.] {Select the objective function to be minimized, denoted in what follows as the primary objective (without loss of generality let us take $J_1$ as the primary objective).}
\item[4.] {For the non-minimized objective $J_2$, generate a vector
              \[ \varepsilon = [ \varepsilon_{1}, \hdots, \varepsilon_{i}, \hdots, \varepsilon_{m} ]\]

              such that $\varepsilon_{1}\leq\underline{J}_2$, $\varepsilon_{m}\geq\overline{J}_2$ and $\varepsilon_{1}<\varepsilon_{2}<\hdots<\varepsilon_{m}$.}
 \item[5.] {Solve the MINLP:
 \begin{equation*}
\min_{x,y}    ~~  J_1 (\dot{z},z,x,y,k)
\end{equation*}
subject to:
\[\varepsilon_{k} \leq J_2(\dot{z},z,x,y,k) < \varepsilon_{k+1} \]
for $k=1,\hdots,m-1$ by means of one of the MINLP solvers introduced in the previous section.}
\item[6.]{Evaluate the solutions obtained and construct the Pareto front with the non dominated optimal ones}.
\end{itemize}

Continuing with the example introduced in the previous section where the goal was to find a circuit with a specific response upon induction, we introduce now an additional design criterium. As mentioned, in the original formulation, the design objective was unique and given by Eq. (\ref{eq:ex1Z}).
Here we consider the protein production cost as an additional objective to minimize, competing with the circuit performance. This criterium has been suggested as  a design principle by several authors \cite{Szekely:2013,Zaslaver:2004}.  The cost of protein production is encoded in an objective function that, taking into account the mass balance equations (\ref{eq:ex1balance}) takes the form:
\begin{equation}\label{eq:cost}
C =  \int_0^T V_j dt
\end{equation}
where $T$ is the final time.

We apply the constraint strategy combined with the MINLP solvers to obtain the Pareto front for different degrees of circuit complexity. First, we set the maximum number of pairs to $M_{max}=2$, and then we increase the maximum number of pairs to evaluate how the Pareto boundary evolves, and how the methodology proposed scales with the systems size. The results obtained are included
in Results and Discussion section.

One interesting application of the methodology presented is to explore new topologies of medium or high order that perform a desired (complex) functionality.
To illustrate this we make use of the same library of components of the previous example, but in this case searching for a circuit topology with the capability to perform adaptation, setting a priori the desired level of complexity.

Adaptation is defined as the ability of the circuit to reset itself after responding to a stimulus. Here, we evaluate the levels of $LacI$ (output) in response to a sustained stimulus of $aTc$ (input). Ma \emph{et al} \cite{Ma:2009} assessed the ability a circuit to adapt after a given stimulus by measuring two functional quantities encoded in two competing objectives related to the sensitivity and the precision of the system's response. On the one hand, in order to maximize adaptation after a given stimulus we need to maximize the circuit's sensitivity:
\begin{equation}\label{eq:sensitivity}
S=O_{peak}-O_{t=0}
\end{equation}
where $O_{peak}$ is the level of the output (in this case $LacI$) at its maximum upon induction and $O_{t=0}$ is the level of the output at the steady state before induction.
On the other hand, in order to maximize adaptation we need to maximize the circuit's precision, i.e. we need to minimize the following function:
\begin{equation}\label{eq:precision}
P=O_{t=T}-O_{t=0}
\end{equation}
where $O_{t=T}$ is the level of the output at steady state reached upon induction. The search for an adaptive circuit can be formulated then as a multiobjective optimization problem where the constraints are imposed by the circuit's dynamics.
In this way, it is possible to elucidate whether is it possible to construct a circuit with capacity from the available set of components.
The maximum and minimum number of allowed components can be adjusted by means of inequality constraints.
The details and results of the corresponding multiobjective optimization problem are included in Results and Discussion section.

\section*{Results and Discussion}
\subsection*{Single objective global optimization  design of a circuit with inducer specific response}

In this section we present the results obtained for the monoobjective problem described in Methods section. Starting from the library of transcripts and promoters indicated we search for the circuit with best circuit performance by maximizing $Z$ defined in Eq. (\ref{eq:ex1Z}). We use the MINLP solvers \texttt{MITS}, \texttt{ACOmi} and \texttt{eSS} with the goal of minimizing $J=-Z$. We solve the optimization problem for increasing levels of complexity, i.e. for an increasing upper bound in the maximum number of pairs ($M_{max}$).
Note that, for a library with $p$ different pairs, the number of possible circuits containing exactly $M$ active pairs is:
\begin{equation}\label{eq:nc}
nC(p,M) = \frac{p!}{(p-M)!M!}.
\end{equation}
According to this formula, the number of combinations $nC$ increases with $p$ and also with the maximum number of pairs $M_{max}$ as illustrated in the Figs. S1 and S2 in the Appendix I. In what follows we do not modify the original library of transcripts and promoters (p=32) and evaluate the performance of the methods for different values of $M_{max}$.

For $M_{max}=2$, the three MINLP solvers, \texttt{MITS}, \texttt{ACOmi} and \texttt{eSS}, reached the same solution, the circuit with active pairs $(Plac_1,tetR)$ and $(Ptet_2,LacI)$. In Fig. 1, we illustrate the best circuit found together with the corresponding superstructure matrix, coefficients of the model and active pairs.
The value of the objective function for the optimal circuit is $J=-0.99998$. This solution coincides with the one obtained by Dasika and Maranas using the outer approximation method \cite{Dasika:2008}.
\begin{figure}[h]
\center
\includegraphics[angle=0,totalheight=6cm]{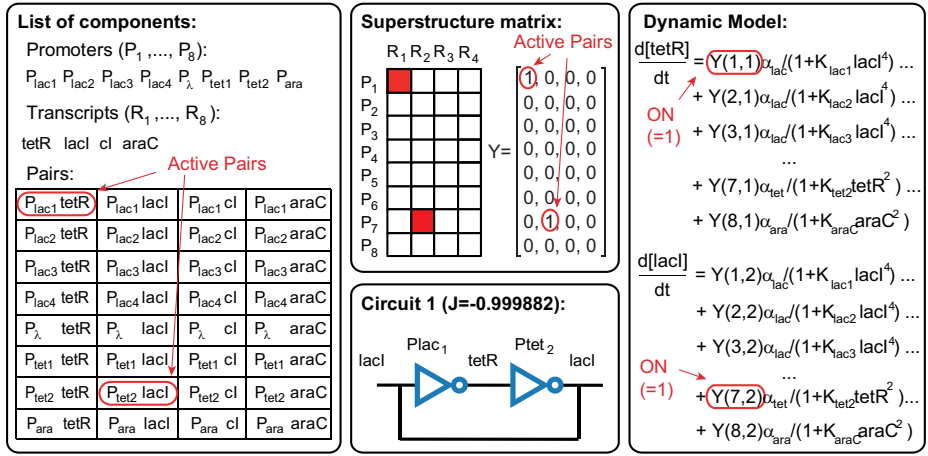}
  \caption{Optimum of the single objective design problem from \cite{Dasika:2008} with the corresponding active pairs, superstructure matrix and dynamic model equations. The full model equations can be found in the Appendix II.}
\end{figure}
The three global MINLP methods achieve the solution in substantially less computation time than the outer approximation method and in particular \texttt{MITS} showed the best performance for this example. Whereas the time reported to find the optimum with the outer approximation method was of 200 minutes in an Intel 3.4 GHz Xeon processor \cite{Dasika:2008}, \texttt{MITS} arrived to the same solution in less than 200 seconds using a slightly slower processor (Intel 2.8 GHz Xeon), thus reducing the computational cost at least by a factor of 60.
To test the algorithm, we have used as starting guess the zero vector $0 \in \mathbb{Z}^{32}$, since the objective function value is very far from the optimum and the constraint is fulfilled. We repeat the analysis starting from different initial guesses fulfilling the constraint and the solver reaches the same solution in similar time. The corresponding convergence curves are illustrated in Fig. S3 from the Appendix I.

Here it is worthy of note that for the monoobjective problem there exist a number of different circuits with similar performance. In Fig. 2 we include four different solutions (circuits 2 to 5) showing very good performance, with values of $J$ below $-0.95$ (note that by definition the minimum value that $J$ can reach is -1). In absence of additional design criteria, and taking into account that different sources of error limit the precision of biocircuit implementations, the selection of the best design among a set of candidates with close objective function values is rather arbitrary.

For $M_{max}=3$, the best solution found is the circuit 6 in Fig. 2, with $J=-0.999996$. Again, \texttt{MITS} showed the best performance, achieving the solution in less than 300 seconds, as it is shown in the convergence curves illustrated in Fig. S4 in the Appendix II.
\begin{figure}[h!]
\center
\includegraphics[angle=0,totalheight=3.9cm]{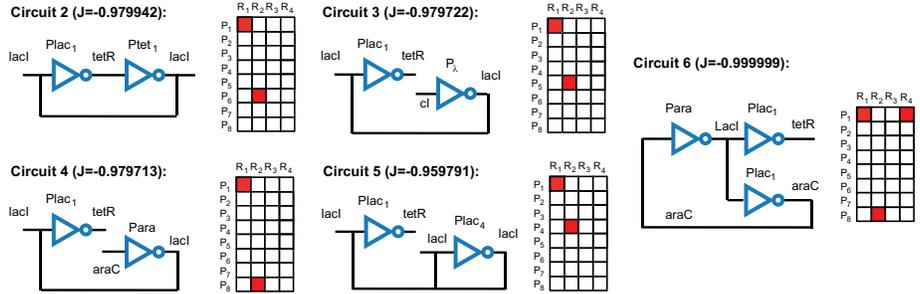}
\caption{Alternative circuit configurations with a maximum of 2 active pairs and similar levels of performance.}
\end{figure}

For $M_{max}=32$, i.e. increasing the maximum level of complexity to 32 pairs (note that this is equivalent to the unconstrained problem), the best solution found is the circuit 7 in Fig. 3 with 14 active pairs. It is important to remark that for increasing levels of complexity the number of solutions with similar values of the objective function (and consequently similar performance) also increases. As an example, we show the circuits 8 to 11 in Fig. 3 with similar level of performance and rather different topologies (for space reasons we depict only the superstructure matrix for all circuits except from 8). Note also that in terms of performance, circuit 7 in Fig. 3 is equivalent to circuit 6 in Fig. 2. This fact leads to arbitrariness when it comes to select the best solution, and suggest the convenience of introducing additional competing criteria in order to provide more realistic design settings.

Regarding solvers performance we observe that, at least for short computation times, the solution found depends on the initial guess (this dependency increases with complexity) and therefore we test every method starting from different initial guesses. Fig. S5 in the Appendix I illustrates the convergence curves of \texttt{MITS} starting from different initial guesses.
Concerning the best circuits found, circuits 7, 8 and 10  in Fig. 3 were obtained by \texttt{MITS} in less than 1500 s, circuit 9 in Fig. 4 was found by \texttt{ACOmi} in less than 3 hours and circuit 11 was found by \texttt{eSS} in less than 300 seconds.

Remarkably, the three methods \texttt{MITS}, \texttt{ACOmi} and Ess provided solutions with objective function values below $J=-0.9999$ in less than 300 seconds, for all the initial guesses tested.
\begin{figure}[h!]
\center
\includegraphics[angle=0,totalheight=6cm]{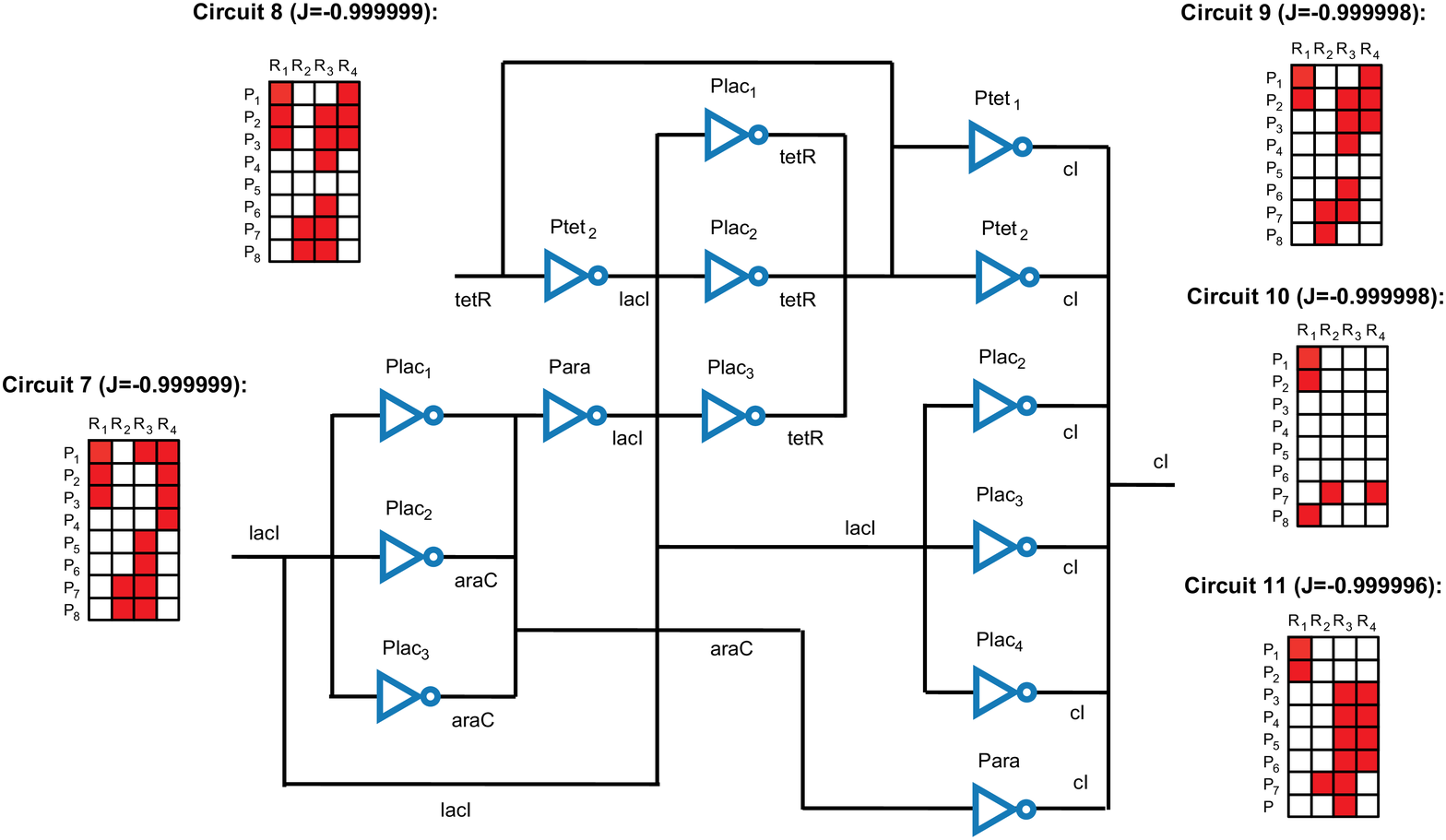}
  \caption{Best circuit found (circuit 7) and  alternative circuits with a maximum of 32 active pairs with similar levels of performance (only superstructure matrices are included except from circuit 8).}
      \end{figure}

\subsection*{Multiobjective global optimization design of a circuit with inducer specific response.}
Next, we introduce the protein production cost as an additional criterium to the design problem. Our primary objective is now the performance function $J_1=-Z$  where $Z$ is given by Eq. (\ref{eq:ex1Z}) and the secondary objective is the cost $J_2=C$, where $C$ has been defined in Eq. (\ref{eq:cost}). The problem is solved for increasing levels of complexity, applying the $\varepsilon$-constraint strategy.

For $M_{max}=2$  we know the solution $y_1^*$ from the previous monoobjective analysis, and the value of the cost at this optimum is $J_2(y_1^*)=2432.3518$. We search now the individual optimum $y_2^*$ for the secondary objective, finding the circuit with active pairs $(Plac_1,LacI)$ and $(Plac_1,tetR)$. Solutions with values of $J_1>0$ are discarded. The value of the cost at the optimum is $J_2(y_2^*)=1129.09$.  Taking into account that the upper and lower bounds for the secondary objective function are
precisely $J_2(y_1^*)=2432.3518$ and $J_2(y_2^*)$, and with a step size of 50, we obtain six non dominated points $P_1 \hdots P_6$ corresponding to six circuits with different topologies.
The Pareto front is illustrated in Fig. 4. The three MINLP solvers have been used in order to compare the results, and an exhaustive search was also implemented, arriving to the same Pareto optimal front. Let us remind that the exhaustive search is possible only for low levels of complexity, since the computation time explodes as the number of maximum pairs increases, as deduced from Eq. \ref{eq:nc}.

\begin{figure}[hb]
\center
\includegraphics[angle=0,totalheight=8cm]{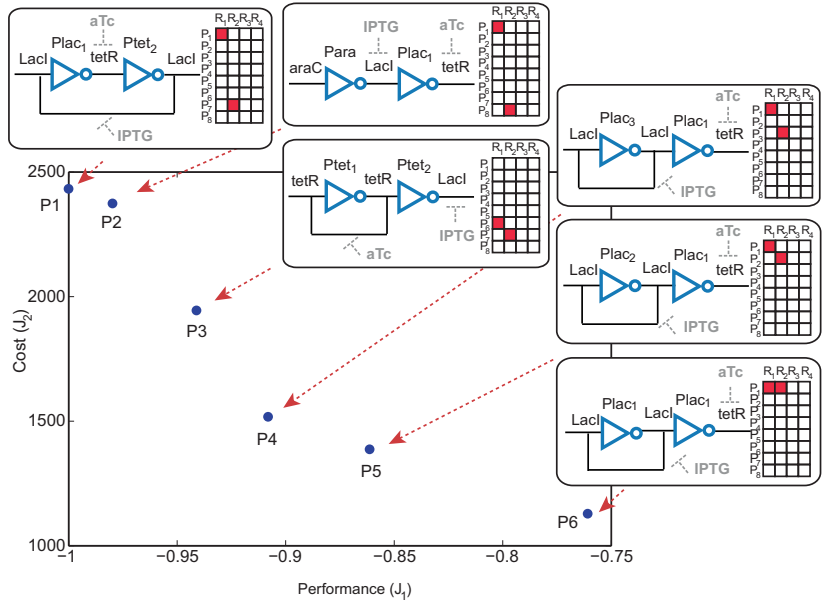}
 \caption{Pareto front for a maximum of 2 active pairs.}
\end{figure}

Following the same strategy, we compute the Pareto front for $M_{max}=3$. The front obtained is shown in Fig. 5, and consists of four different points, labeled $Q_1 \hdots Q_4$ (note that $Q_2 = P_1$).

It is of relevance that the solution $Q_4$ is significantly better in terms of cost than any other and at the same time it shows a very good performance $(J_1 < -0.95)$. The multiobjective formulation allowed in this case to find a non intuitive topology which is a very good candidate for a successful laboratory implementation. It can be deduced also from Fig. 5
how an small increase in complexity from $M_{max}=2$ to $M_{max}=3$ leaded to significant improvement in the Pareto front, where $Q_3$ and $Q_4$ are non dominated by any of the circuits with two active pairs ($P_1 \hdots P_6$).
\begin{figure}[h!]
\center
\includegraphics[angle=0,totalheight=7.5cm]{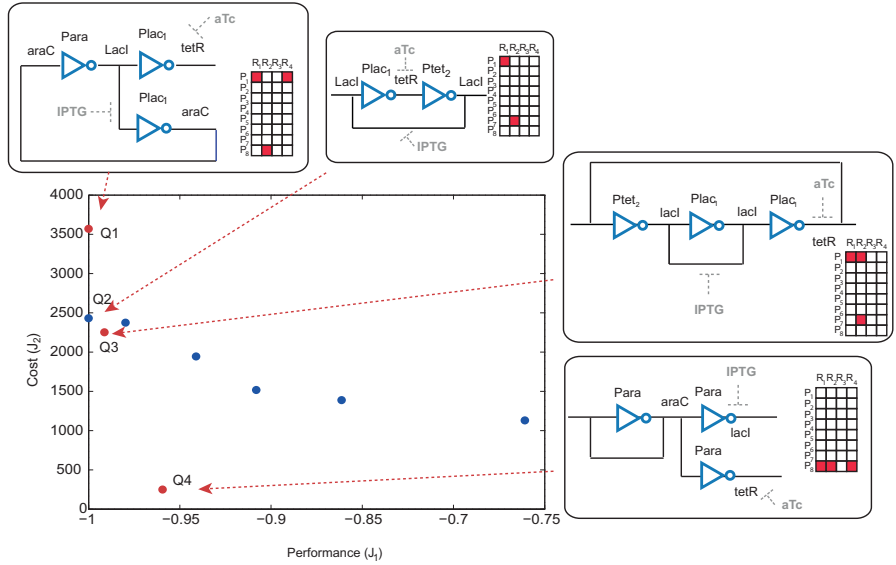}
  \caption{Pareto front for  a maximum of 3 active pairs.}
\end{figure}

Finally, we compute the Pareto front for $M_{max}=32$. The circuit $Q_1$ (circuit 6 in Fig. 2) obtained for $M_{max}=3$ is also the best solution found for the unconstrained problem (together with the circuit 7 in Fig. 3). By constraining also the minimum level of complexity by setting  $M_{max}>3$ we obtain the set of non-dominated solutions depicted in
Fig. 6, together with the corresponding superstructure matrices. In this figure it can be seen that the multiobjective strategy employed allowed us to find points in non-convex regions of the Pareto front, as it is the case of the circuit $R_5$.
\begin{figure}[h!]
\center
\includegraphics[angle=0,totalheight=8cm]{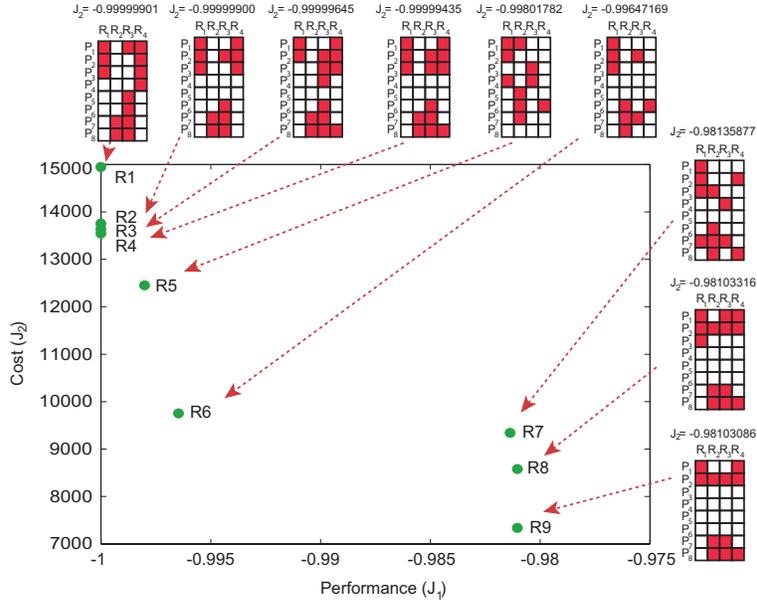}
  \caption{Pareto front for a minimum of 3 active pairs and a maximum of 32 active pairs.}
\end{figure}
\subsection*{Adaptive biocircuit with predefined complexity}

Now, starting from the same library of components of the previous example we search for a circuit configuration with the ability for adaptation.  We assume that one of the parameters can be manipulated, in this case a kinetic parameter related to the $P_{tet}$ promoter $\alpha_{tet}$ (see Appendix II).
As indicated in  Methods section  the adaptive capacity of the circuit is evaluated by the levels of the output protein $LacI$ in response to a sustained stimulus of $aTc$, in particular by the levels at its maximum upon induction $O_{peak}$,  at the steady state before induction $O_{t=0}$ and at the steady state upon induction $O_{t=T}$. Two competing objective functions are considered, the circuit's sensitivity defined by Eq. (\ref{eq:sensitivity}) and the circuit's precision measured through the formula in Eq. (\ref{eq:precision}).
The multiobjective MINLP problem with 32 integer and 1 real decision variables is solved with the $\varepsilon$-constraint strategy proposed,  maximizing as a primary objective the sensitivity, i.e. minimizing $-S$ with $S$  defined in Eq. (\ref{eq:sensitivity}) and setting the precision as a constraint. In Fig. 7 A, we depict one of the solutions of the Pareto front, where $P<20$ with $P$ defined in Eq. (\ref{eq:precision}). As it is shown in Fig. 7 B, the circuit is able to adapt upon a sustained stimulus of $aTc$. The optimal value for the kinetic constant is also indicated.
\begin{figure}[h!]
\center
\includegraphics[angle=0,totalheight=3.8cm]{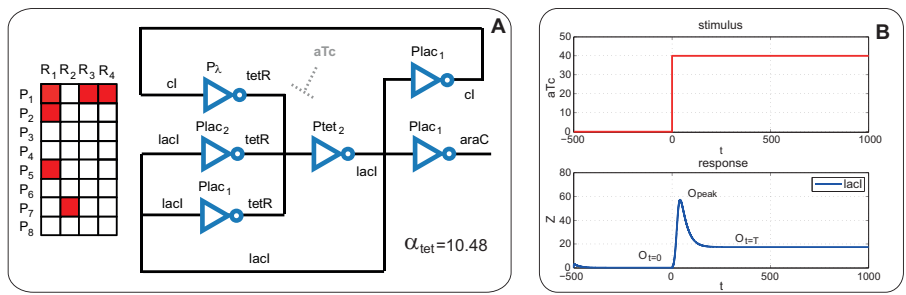}
  \caption{Adaptive circuit found by multiobjective optimization. A) Topology and superstructure matrix. B) Circuit's response upon aTc stimulus.}
\end{figure}

\section*{Conclusions}
In this work we have introduced a multiobjective formulation for the design of biocircuits. The presence of more than one competing objective provides  more realistic design settings where the solution is not unique and every solution represents trade-off between different criteria.

The multiobjective optimization in the context of genetic circuit design posed a number of challenges mainly due to the inherent nonlinear nature of the gene circuit's dynamics and the large search spaces involved combining the presence of integer and real variables, which makes the expected Pareto front discrete and possibly non-convex.

In order to overcome these difficulties we made use of global optimization algorithms, showing their efficiency for the MINLP problem resultant of the monoobjective formulation of the design.
Then, we provided a multiobjective optimization framework for the design of biocircuits that combines the efficiency of the global MINLP solvers with the capacity to handle multiple design criteria.

Looking for further extensions the method presented is quite general, accommodating to any ODE based modeling framework such that the circuit's model structure is obtained from the starting list of parts by giving values to a set of binary variables.

The advantages of this multiobjective formulation were shown through the design of a biocircuit with specific response upon induction. Due to the efficiency of the global solvers it was possible to obtain in reasonable times the Pareto fronts for different levels of complexity including circuits belonging to non-convex regions of the optimal set of solutions. The capacity to handle circuits with higher number of regulatory regions implies more opportunities for parameter tuning.

Through an illustrative example, we have demonstrated how using this framework we can obtain non intuitive designs to perform a desired functionality setting up a priori the desired level of
complexity. This can be useful in future contributions to explore and identify different design principles for synthetic gene circuits.

\bibliographystyle{unsrtnat} 
\bibliography{Otero-Muras-and-Banga}      

\newpage
\appendix{\textbf{Appendix I}}
This appendix contains the generic equations for the model of the genetic circuit used as a working example through the main text. The dynamic model equations and the kinetic parameters have been obtained from Ref. \cite{Dasika:2008}.

The states of the generic model are the concentrations of all the species involved $lacI$, $LacIIPTG$, $tetR$, $tetRaTc$, $cI$, $araC$ and the generic model consists of 6 ordinary differential equations representing the mass balances for the system. The coefficients $Y(i,j)$ for $i=1,\hdots,8$ and $j=1,\hdots,4$ are binary variables that take value 1 if the corresponding pair is active and 0 if it is inactive for a given circuit.
\begin{align*}
dlacI = &+ Y(1,2)\alpha_{lac}/(1+K_{lac1}lacI^4) \\
        &+ Y(2,2)\alpha_{lac}/(1+K_{lac2}lacI^4) \\
        &+ Y(3,2)\alpha_{lac}/(1+K_{lac3}lacI^4) \\
        &+ Y(4,2)\alpha_{lac}/(1+K_{lac4}lacI^4) \\
        &+ Y(5,2)\alpha_{\lambda}/(1+K_{lambda} cI^2) \\
        &+ Y(6,2)\alpha_{tet}/(1+K_{tet1}tetR^2) \\
        &+ Y(7,2)\alpha_{tet}/(1+K_{tet2}tetR^2) \\
        &+ Y(8,2)\alpha_{ara}/(1+K_{araC}araC^2) \\
        &- Kf \cdot lacI \cdot IPTG \\
        &+ Kb \cdot lacIIPTG \\
        &- K_{deglacI} \cdot lacI \\
 dlacIIPTG &= + Kf\cdot lacI \cdot IPTG \\
        &- Kb\cdot lacIIPTG  \\
        &- K_{degcpx}\cdot lacIIPTG \\
dtetR = &+ Y(1,1) \alpha_{lac}/(1+K_{lac1}lacI^4) \\
        &+ Y(2,1) \alpha_{lac}/(1+K_{lac2}lacI^4) \\
        &+ Y(3,1) \alpha_{lac}/(1+K_{lac3}lacI^4) \\
        &+ Y(4,1) \alpha_{lac}/(1+K_{lac4}lacI^4) \\
        &+ Y(5,1) \alpha_{\lambda}/(1+K_{lambda}  \cdot cI^2) \\
        &+ Y(6,1) \alpha_{tet}/(1+K_{tet1}tetR^2) \\
        &+ Y(7,1) \alpha_{tet}/(1+K_{tet2}tetR^2) \\
        &+ Y(8,1) \alpha_{tet}/(1+K_{araC}araC^2) \\
       & - Kf tetR\cdot aTc \\
        &+ Kb\ tetRaTc \\
       & - K_{degtetR}tetR \\
dtetRaTc = &+ Kf\cdot tetR\cdot aTc \\
        &- Kb\cdot tetRaTc \\
        &- K_{degcpx}\cdot tetRaTc \\
dcI = &+ Y(1,3) \alpha_{lac}/(1+K_{lac1}lacI^4) \\
        &+ Y(2,3)\alpha_{lac}/(1+K_{lac2}lacI^4) \\
        &+ Y(3,3)\alpha_{lac}/(1+K_{lac3}lacI^4) \\
        &+ Y(4,3)\alpha_{lac}/(1+K_{lac4}lacI^4) \\
        &+ Y(5,3) \alpha_{\lambda}/(1+K_{lambda}  \cdot cI^2) \\
        &+ Y(6,3)\alpha_{tet}/(1+K_{tet1}tetR^2) \\
        &+ Y(7,3)\alpha_{tet}/(1+K_{tet2}tetR^2) \\
        &+ Y(8,3)\alpha_{tet}/(1+K_{araC}araC^2) \\
        &- K_{degcI} cI\\
daraC = &+ Y(1,4)\alpha_{lac}/(1+K_{lac1}lacI^4) \\
        &+ Y(2,4)\alpha_{lac}/(1+K_{lac2}lacI^4) \\
        &+ Y(3,4)\alpha_{lac}/(1+K_{lac3}lacI^4) \\
        &+ Y(4,4)\alpha_{lac}/(1+K_{lac4}lacI^4) \\
        &+ Y(5,4)\alpha_{\lambda}/(1+K_{lambda}  \cdot cI^2) \\
        &+ Y(6,4)\alpha_{tet}/(1+K_{tet1}tetR^2) \\
        &+ Y(7,4)\alpha_{tet}/(1+K_{tet2}tetR^2) \\
        &+ Y(8,4)\alpha_{tet}/(1+K_{araC}araC^2) \\
        &- K_{degaraC}araC.
\end{align*}

The values of the parameters are included in the following Table S1.

\setcounter{figure}{0}
\setcounter{table}{0}
\renewcommand{\thetable}{S\arabic{table}}
\renewcommand{\thefigure}{S\arabic{figure}}
\begin{table}[h!]
\center
\caption{Parameters}
\begin{tabular}{lcl}
\hline
Parameter & value \\
\hline
$\alpha_{lac}$ &  .215 \\
$\alpha_{tet}$ & 1.215 \\
$\alpha_{\lambda}$ & 2.92 \\
$\alpha_{ara}$ & 1.215 \\
$K_{\lambda}$ & 0.33 \\
$K_{tet1}$ & 0.014 \\
$K_{tet2}$ & 1.4 \\
$K_{lac1}$ & 10 \\
$K_{lac2}$ & 0.01 \\
$K_{lac3}$ & 0.001 \\
$K_{lac4}$ & 0.00001\\
$K_{araC}$ & 2.5 \\
$K_{deglacI}$ & 0.0346 \\
$K_{degtetR}$ & 0.0346 \\
$K_{degcI}$ & 0.0693 \\
$K_{degaraC}$ & 0.0115 \\
$K_{degcpx}$ & 0.0693 \\
$K_{f}$ & 0.05 \\
$K_{b}$ & 0.1 \\
\hline
\end{tabular}
\end{table}

\newpage

\appendix{\textbf{Appendix II}}
This appendix contains additional figures S1, S2, S3, S4 and S5.
\beginsupplement
\begin{figure}[h!]
\begin{center}
\includegraphics[angle=0,totalheight=6cm]{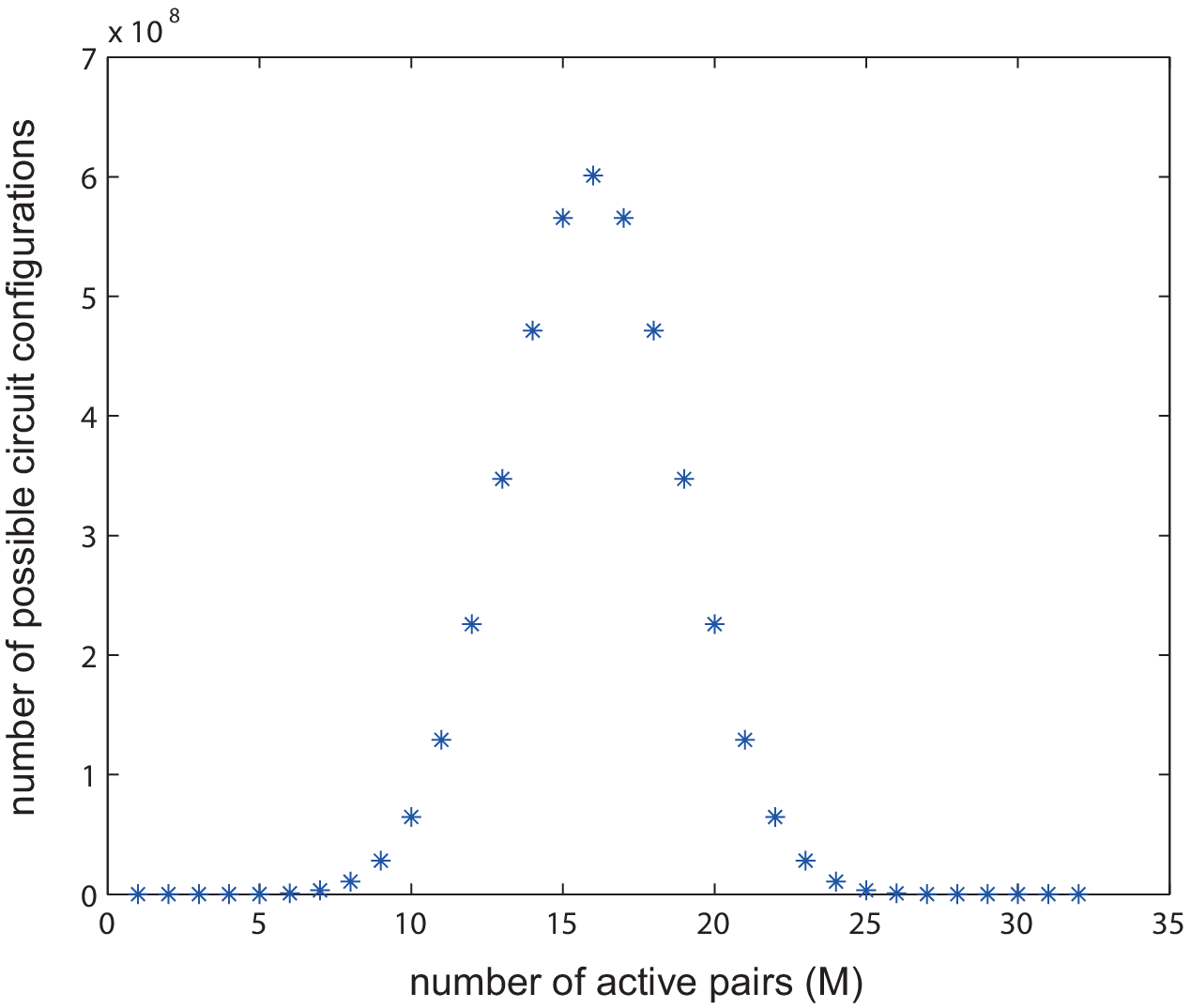}
\caption{Number of possible circuit configurations as a function of the number of active pairs.}  
\label{figS1}
\end{center}
\end{figure}
\begin{figure}[h!]
\begin{center}
\includegraphics[angle=0,totalheight=6cm]{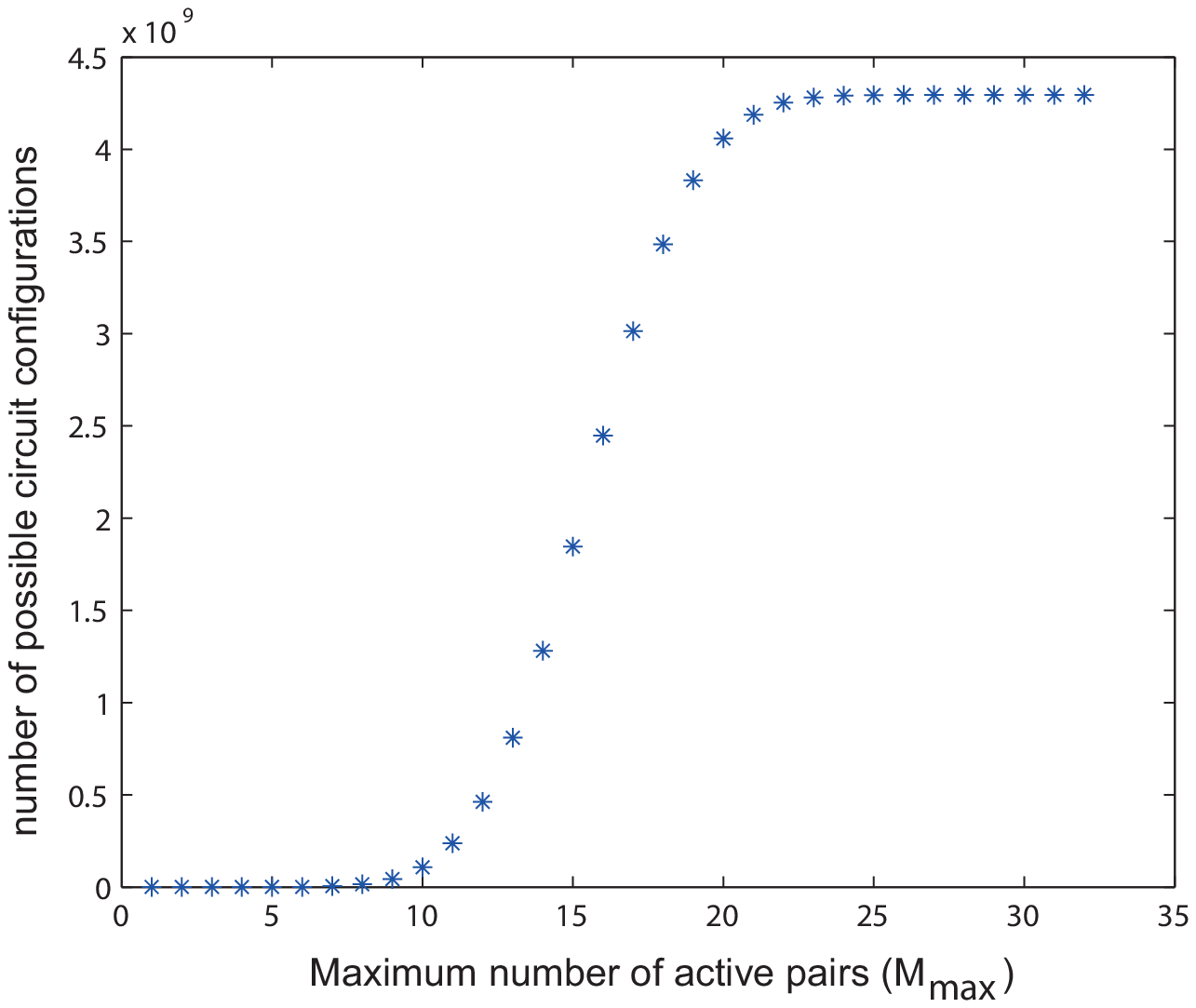}
\caption{Number of possible circuit configurations for increasing $M_{max}$.}  
\label{figS1}
\end{center}
\end{figure}
\begin{figure}[h!]
\begin{center}
\includegraphics[angle=0,totalheight=6cm]{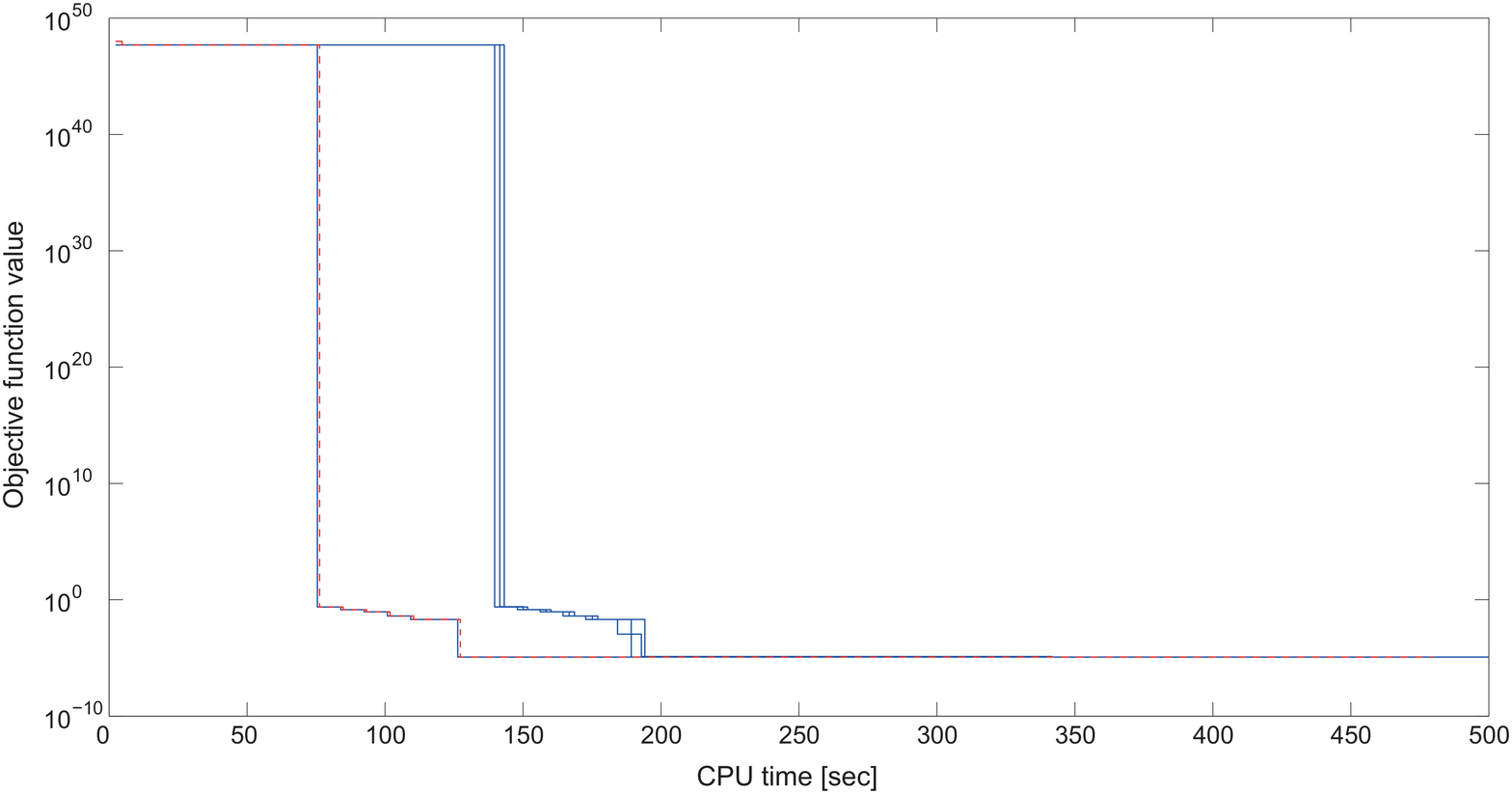}
\caption{MITS Convergence curves for different initial guesses with $M_{max}=2$. The initial guess for the dashed convergence curve is the circuit with no active pairs. To make possible
a log scale representation the objective function has been shifted to the positive orthant.}  
\label{figS1}
\end{center}
\end{figure}
\begin{figure}[h!]
\begin{center}
\includegraphics[angle=0,totalheight=6cm]{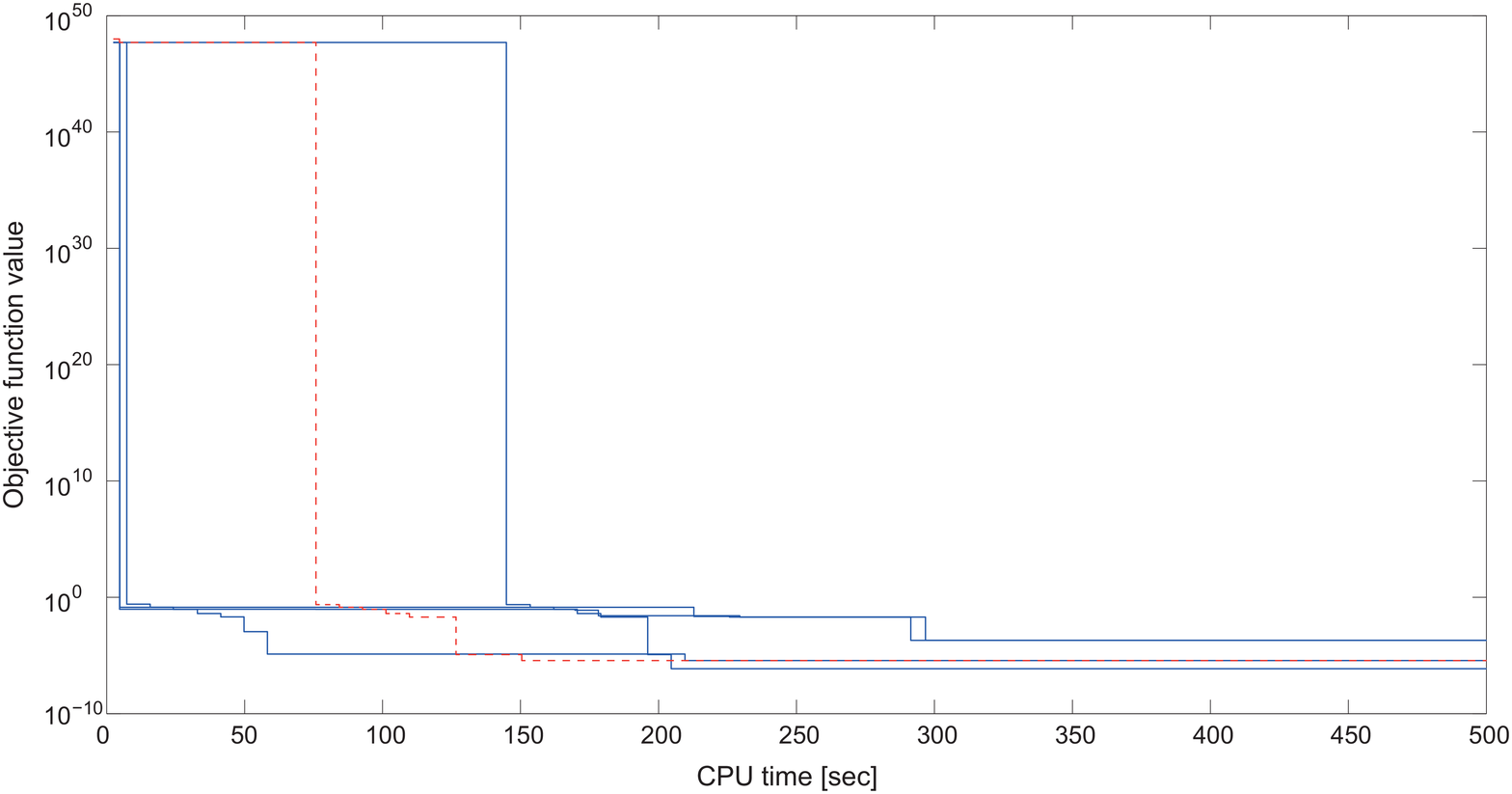}
\caption{MITS Convergence curves for different initial guesses with $M_{max}=3$. The initial guess for the dashed convergence curve is the circuit with no active pairs. To make possible
a log scale representation the objective function has been shifted to the positive orthant.}  
\label{figS2}
\end{center}
\end{figure}
\begin{figure}[h!]
\begin{center}
\includegraphics[angle=0,totalheight=6cm]{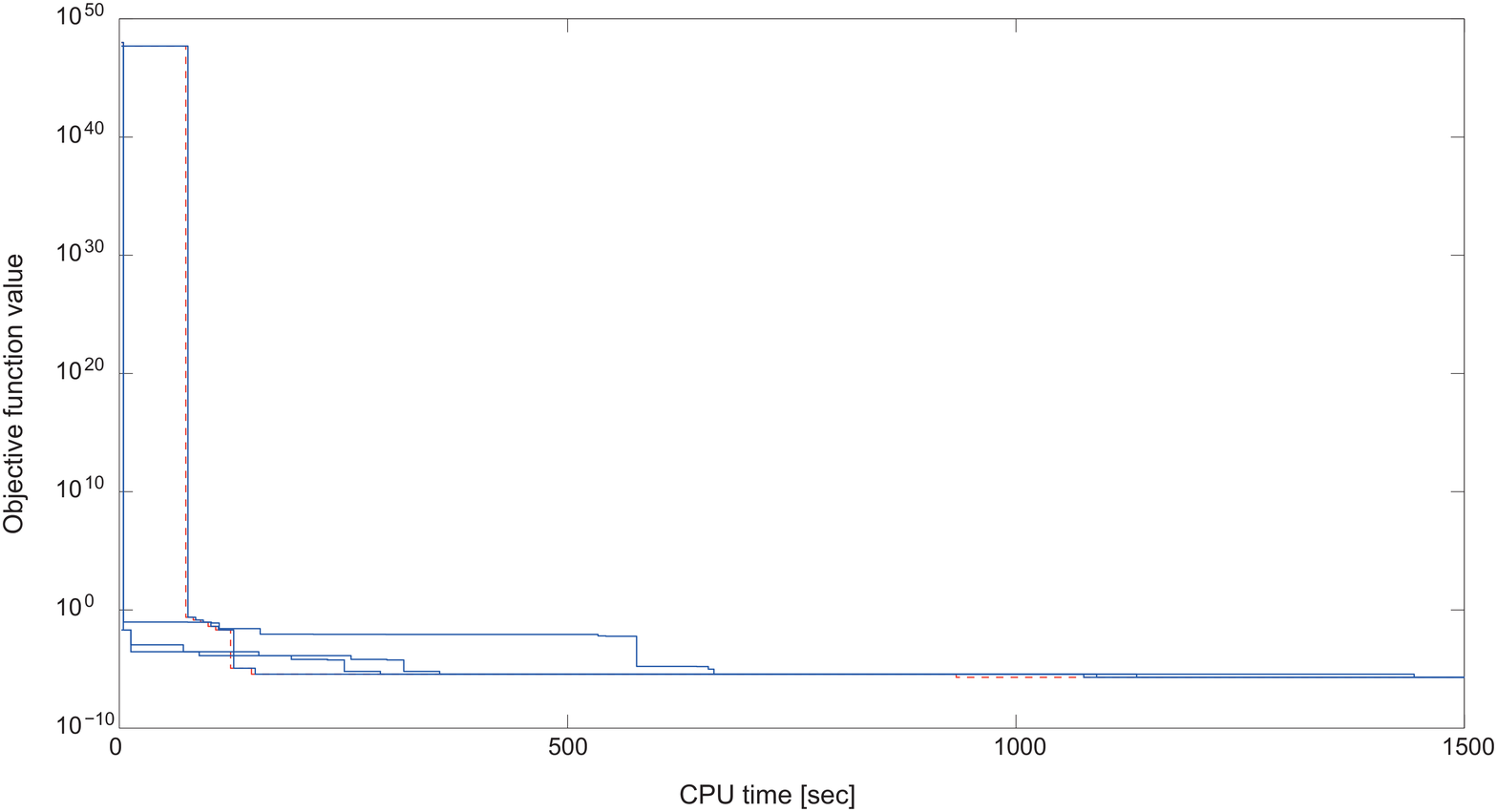}
\caption{MITS Convergence curves for different initial guesses with $M_{max}=32$. The initial guess for the dashed convergence curve is the circuit with no active pairs. To make possible
a log scale representation the objective function has been shifted to the positive orthant.}  
\label{figS2}
\end{center}
\end{figure}

\end{document}